\documentclass[10pt]{article}
\usepackage{amsthm}
\usepackage{amsfonts}
\usepackage{amssymb,latexsym}
\usepackage{amsbsy}
\usepackage{amsxtra}
\usepackage{hyperref}
\usepackage{amsmath}
\usepackage{mathabx}
\usepackage{mathrsfs}
\usepackage{graphicx}
\usepackage{bbm}

\theoremstyle{plain}
\newtheorem{theorem}{Theorem}
\newtheorem{lemma}{Lemma}
\newtheorem{proposition}{Proposition}

\begin{document}

\title{$(L^{r}, L^{s})$ Resolvent Estimate for the Sphere off the Line $\frac{1}{r}-\frac{1}{s}=\frac{2}{n}$}
 \author{Tianyi Ren}
 
\maketitle
\begin{abstract}
We extend the resolvent estimate on the sphere to exponents off the line $\frac{1}{r}-\frac{1}{s}=\frac{2}{n}$. Since the condition $\frac{1}{r}-\frac{1}{s}=\frac{2}{n}$ on the exponents is necessary for a uniform bound, one cannot expect estimates off this line to be uniform still. The essential ingredient in our proof is an $(L^{r}, L^{s})$ norm estimate on the operator $H_{k}$ that projects onto the space of spherical harmonics of degree $k$. In showing this estimate, we apply an interpolation technique first introduced by Bourgain \cite{Bourgain}. The rest of our proof parallels that in Huang-Sogge \cite{HS}.

\end{abstract} 

\maketitle
\section{Introduction}

Since Kenig-Ruiz-Sogge's classical result \cite{KRS} on resolvent estimate in the Euclidean space, there has been a lot of interest and endeavor in extending this work to manifolds. Their classical result says that in $\mathbb{R}^{n}$ where $n \geqslant 3$, the uniform estimate \[
\|u\|_{L^{s}(\mathbb{R}^{n})} \leqslant C\|(\Delta+\zeta)u\|_{L^{r}(\mathbb{R}^{n})}, \quad u \in H^{2, r}(\mathbb{R}^{n}), \quad \zeta \in \mathbb{C} \]
holds for pairs of exponents $r$ and $s$ satisfying \begin{equation} \label{29}
\frac{1}{r}-\frac{1}{s}=\frac{2}{n}, \quad \frac{2n}{n+3}<r<\frac{2n}{n+1}.
\end{equation}
Here uniformity means that the constant $C$ is independent of the function $u$ and in particular, of $\zeta \in \mathbb{C}$. 

In 2011, Dos Santos Ferreira, Kenig and Salo \cite{FKS} showed that on a compact manifold $M$, the uniform inequality \[
\|u\|_{L^{\frac{2n}{n-2}}(M)} \leqslant C\|(\Delta_{g}+\zeta)u\|_{L^{\frac{2n}{n+2}}(M)} \]
is valid for all $\zeta \in \mathbb{C}$ such that $\mathrm{Im}\sqrt{\zeta} \geqslant \delta$ for a fixed positive $\delta$. The two exponents in their result, duals of each other, correspond to the midpoint of the line segment \eqref{29} in Kenig-Ruiz-Sogge \cite{KRS}. On a manifold of course, the Laplacian becomes the Laplace-Beltrami operator. Pictorially, the permissible $\zeta \in \mathbb{C}$ constitute the region in the complex plane outside a small ball around the origin and a parabola. This naturally poses the question of whether we are able to enlarge the region for $\zeta$, while still obtaining a uniform estimate. 

In 2013, Bourgain, Shao, Sogge and Yao \cite{BSSY} showed that the region in Dos Santos Ferreira-Kenig-Salo \cite{FKS} is optimal for Zoll manifolds, hence in particular for spheres. However, for the torus and manifolds of nonpositive curvature, they expand the region remarkably. Later, Shao and Yao \cite{SY} proved the resolvent estimate on compact manifolds for exponents that do not necessarily lie on the line of duality. More explicitly, in their theorem, the exponents $r$ and $s$ need only satisfy \[
\frac{1}{r}-\frac{1}{s}=\frac{2}{n}, \quad \frac{2(n+1)}{n+3} \leqslant r \leqslant \frac{2(n+1)}{n-1}. \]
This leads us to the question of whether it is possible to extend the permissible exponents to a greater part of the line segment \eqref{29}, besides the question of finding the sharp region for $\zeta$ on certain types of manifolds. 

In 2014, Huang and Sogge \cite{HS} proved the resolvent estimate on spaces of constant positive and negative curvature for exponents lying on the full line segment \eqref{29}. The region for $\zeta$ in the case of constant positive curvature is the same as that in Dos Santos Ferreira-Kenig-Salo \cite{FKS} since it was shown to be sharp, as just mentioned. The region in the case of constant negative curvature is the whole complex plane when the dimension is $3$, and the complex plane with a neighborhood of the origin excluded in higher dimensions.

There has also been efforts to prove the estimate for exponents that are not even on the line $\frac{1}{r}-\frac{1}{s}=\frac{2}{n}$. By an easy dilation argument, as Kenig-Ruiz-Sogge \cite{KRS} pointed out, the condition $\frac{1}{r}-\frac{1}{s}=\frac{2}{n}$ is necessary in order to obtain a uniform bound independent of $\zeta \in \mathbb{C}$. Therefore, one cannot expect to attain a uniform inequality in this case. However, a bound that depends on $\zeta$ would still be of great interest, especially if it is a negative power of $\zeta$, as this bound will tend to $0$ when $|\zeta|$ goes to infinity. A recent paper by Frank and Schimmer \cite{FS} contributed just to that. They provided the following estimate on compact manifolds: \[
\|u\|_{L^{\frac{2(n+1)}{n-1}}(\mathbb{R}^{n})} \leqslant C|\zeta|^{-\frac{1}{n+1}}\|(\Delta+\zeta)u\|_{L^{\frac{2(n+1)}{n+3}}(\mathbb{R}^{n})}, \]
where as in Dos Santos Ferreira-Kenig-Salo \cite{FKS}, $\mathrm{Im}\sqrt{\zeta} \geqslant \delta$ for a fixed $\delta$, and $C$ of course, is independent of $\zeta$. The exponents, duals of each other, lie on the line $\frac{1}{r}-\frac{1}{s}=\frac{2}{n+1}$.

Frank and Schimmer's result \cite{FS} is for general compact manifolds and concerns the exponents on the line of duality. We in this article, treat the spheres only, but consider more general exponents, those not necessarily on the line $\frac{1}{r}-\frac{1}{s}=\frac{2}{n}$, and not necessarily duals of each other. Our main result is the following

\begin{theorem}
Let $S^{n}$ be the $n$ dimensional sphere where $n \geqslant 3$. If the exponents $r$ and $s$ satisfy \[
\frac{1}{r}-\frac{1}{s}=\sigma, \quad \frac{2}{n+1} \leqslant \sigma \leqslant \frac{2}{n}, \]
\[\frac{2n}{n-1+2n\sigma} < r < \frac{2n}{n+1}, \]
then we have the inequality \begin{equation}
\|u\|_{L^{s}(S^{n}, \mathrm{d}V_{S^{n}})} \leqslant C\lambda^{n\sigma-2}\Big\| \Big( \big( \Delta_{S^{n}}-(\frac{n-1}{2})^{2}\big)+\zeta\Big)u\Big\|_{L^{r}(S^{n}, \mathrm{d}V_{S^{n}})},
\end{equation}
where $\zeta=(\lambda+i\mu)^{2}$, and $\lambda \geqslant 1$, $|\mu| \geqslant 1$.
\end{theorem}

Notice that when $\frac{2}{n+1} \leqslant \sigma < \frac{2}{n}$, the power $n\sigma-2$ on $\lambda$ in the estimate is negative. On the $(\frac{1}{r}, \frac{1}{s})$ plane, the exponents in the theorem constitute the line segments connecting a point whose horizontal axis is $\frac{n+1}{2n}$ with its dual, the point whose vertical axis is $\frac{n-1}{2n}$; see Figure \ref{fig1} below. We, as in Huang-Sogge \cite{HS}, shift the Laplacian to $\Delta_{S^{n}}-(\frac{n-1}{2})^{2}$, because we can then take the square root of minus this shifted Laplacian, and the eigenvalues of the square root are $k+\frac{n-1}{2}$, $k=0, 1, 2, \cdots$. That we are not able to prove the result for $\zeta$ in the optimal region \[
\mathcal{R}=\{\zeta\in\mathbb{C}: \mathrm{Re}\zeta \leqslant (\mathrm{Im}\zeta)^{2}\} \]
is because the case for small $\lambda$ or $|\mu|$ were resolved by Sobolev Embedding Theorem in Huang-Sogge \cite{HS}, but when the exponents are off the line $\frac{1}{r}-\frac{1}{s}=\frac{2}{n}$, the theorem no longer applies.

\begin{figure}
  \centering
    \includegraphics[height=6.5cm]{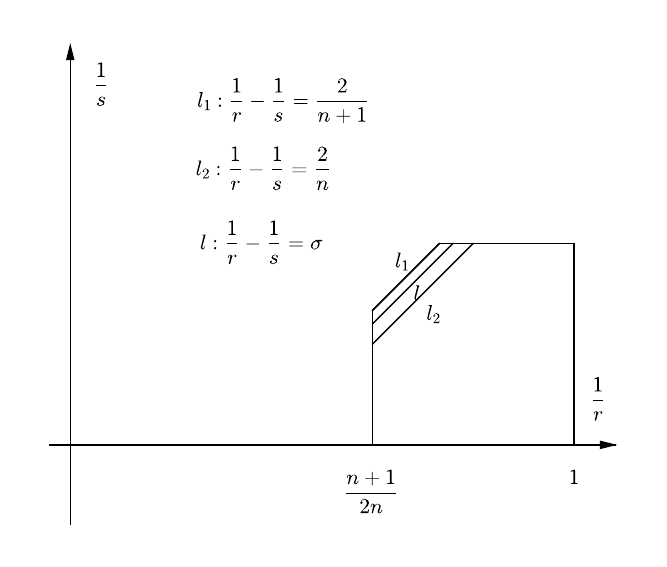}
\caption{Exponents for the resolvent estimate in Theorem 1}
  \label{fig1}
\end{figure}

Let $H_{k}$ denote the projection operator onto the space of spherical harmonics of degree $k$, $k=0, 1, 2, \cdots$. These are the eigenspaces of the square root of the shifted Laplacian, $\sqrt{-\Delta_{S^{n}}+(\frac{n-1}{2})^{2}}$, with eigenvalues $k+\frac{n-1}{2}$. The essential ingredient in proving Theorem 1 is the estimate on the norm of $H_{k}$ stated below.

\begin{proposition}
Let $n \geqslant 3$. We have \begin{equation}
\|H_{k}\|_{L^{r}(S^{n})\rightarrow L^{s}(S^{n})} \leqslant Ck^{n\sigma-1},
\end{equation}
if the exponents $r$ and $s$ are as in Theorem 1, i.e.
\[\frac{1}{r}-\frac{1}{s}=\sigma, \quad \frac{2}{n+1} \leqslant \sigma \leqslant \frac{2}{n}, \]
\[\frac{2n}{n-1+2n\sigma} < r < \frac{2n}{n+1}. \]
\end{proposition}

Finding the norm of $H_{k}$ as an operator from $L^{r}(S^{n})$ to $L^{s}(S^{n})$ when $\frac{1}{r}-\frac{1}{s}=\sigma$ is interesting in its own right. For example, when $\sigma=1$, i.e. $r=1$ and $s=\infty$, then $\|H_{k}\|_{L^{1}(S^{n})\rightarrow L^{\infty}(S^{n})}$ is bounded by $Ck^{n-1}$, because we have the well-known estimate on the $L^{\infty}$ norm of the kernel $H_{k}(x, y)$ of $H_{k}$ (see Sogge \cite{Sogge}, $\S$ 3.4): \[
|H_{k}(x, y)| \leqslant Ck^{n-1}. \]
And of course, if $r=2$, $s=2$, then $\|H_{k}\|_{L^{2}(S^{n})\rightarrow L^{2}(S^{n})}=1$. Determining a bound for $\|H_{k}\|_{L^{r}(S^{n})\rightarrow L^{s}(S^{n})}$ when $\frac{1}{r}-\frac{1}{s}=\sigma$ in the form of a function of $\sigma$, e.g. $k$ raised to a power that is a function of $\sigma$, seems then like a pretty interesting problem. The above proposition is a try in this direction.

To prove the proposition, we utilize an interpolation technique first introduced by Bourgain \cite{Bourgain} when he was proving a bound for the spherical maximal operator. See also Bak-Seeger \cite{BS}, and Carbery-Seeger-Waigner-Wright \cite{CSWW} which contains an abstract version of this result. Bourgain's interpolation technique will be used to prove restricted weak type inequality at the two endpoints of a fixed line segment in the proposition. Then, real interpolation yields the desired estimate on the segment between the endpoints. We already followed this path in \cite{RXZ}, to prove an endpoint version of Stein-Tomas Fourier restriction theorem.

\maketitle
\section{Proof of Proposition 1}

Fix one line segment for the exponents: \begin{equation} \label{30}
\frac{1}{r}-\frac{1}{s}=\sigma, \quad \frac{2n}{n-1+2n\sigma} < r < \frac{2n}{n+1}.
\end{equation}
As just mentioned, we apply Bourgain's interpolation technique to prove restricted weak type estimate at the endpoints. The technique says

\begin{lemma}
If an operator $T$ between function spaces is the sum of the operators $T_{j}$: $T=\displaystyle\sum_{j=1}^{\infty}T_{j}$, and if  each $T_{j}$ satisfies the estimates \begin{equation} \label{1}
\|T_{j}\|_{L^{p_{1}}\to L^{q_{1}}} \leqslant M_{1}2^{\beta_{1}j},
\end{equation}
\begin{equation} \label{2}
\|T_{j}\|_{L^{p_{2}}\to L^{q_{2}}} \leqslant M_{2}2^{-\beta_{2}j},
\end{equation}
for some constants $M_{1}>0$, $M_{2}>0$ and $\beta_{1}>0$, $\beta_{2}>0$, then $T$, the sum of the $T_{j}$, enjoys restricted weak type inequality between the intermediate spaces: \begin{equation} \label{3}
\|T_{j}\|_{L^{p, 1}\to L^{q, \infty}} \leqslant CM_{1}^{\theta}M_{2}^{1-\theta}, \quad where
\end{equation}
\[\theta=\frac{\beta_{2}}{\beta_{1}+\beta_{2}}, \]
\[\frac{1}{p}=\frac{\theta}{p_{1}}+\frac{1-\theta}{p_{2}}, \quad \frac{1}{q}=\frac{\theta}{q_{1}}+\frac{1-\theta}{q_{2}}, \]
and $C$ depends only on $\beta_{1}$ and $\beta_{2}$.
\end{lemma}

Frank and Schimmer's idea in their proof \cite{FS} of the bound for the Hada-mard parametrix actually directs us to a proof of Bourgain's interpolation technique in this special setting of Lebesgue spaces. We record their idea here.

\paragraph{Proof of Lemma 1}
By the property of Lorentz spaces (see Stein and Weiss \cite{SW} $\S$ 5.3, Theorem 3.13 for instance), it suffices to prove \eqref{3} for characteristic functions, i.e. we need only show the inequality \begin{equation} \label{6}
\Big( \mathrm{sup}_{\lambda>0}  \lambda^{q}\mu\big(\{x: |T\mathbbm{1}_{E}|>\lambda\}\big)\Big)^{\frac{1}{q}} \leqslant CM_{1}^{\theta}M_{2}^{1-\theta}\mu(E)^{\frac{1}{p}},
\end{equation}
where $E$ is a measurable set of finite measure, and $\mu(\cdot)$ indicates the measure of a set. For convenience, we denote \[
A=\{x: |T\mathbbm{1}_{E}|>\lambda\}. \]
In what follows, it is necessary to ensure that $A$ has finite measure, but this is in fact a simple consequence of the second estimate \eqref{2} for $T_{j}$ in the Lemma. To see why, first note that by summing the $T_{j}$, we conclude that $T$ is a bounded operator from $L^{p_{2}}$ to $L^{q_{2}}$. This fact together with Tchebyshev's inequality then immediately yields \[
\lambda^{q_{2}}\mu(A) \leqslant \int \big(|T\mathbbm{1}_{E}|\big)^{q_{2}} \mathrm{d}\mu \leqslant \mu(E)^{\frac{q_{2}}{p_{2}}} < \infty. \]

Still by Tchebyshev's inequality, we have \begin{equation}
\mu(A) \leqslant \frac{1}{\lambda}\int_{A} |T\mathbbm{1}_{E}|\mathrm{d}\mu.
\end{equation}
The trick now is to split the sum $T=\displaystyle\sum_{j=1}^{\infty} T_{j}$ into two parts, a finite one and the rest infinite one:
\[ T^{(1)}=\displaystyle\sum_{j=1}^{\rho} T_{j}, \quad T^{(2)}=\displaystyle\sum_{j=\rho+1}^{\infty} T_{j}, \]
where $\rho$ is to be specified later. Then by Holder's inequality, \begin{equation} \label{4} \begin{aligned}
\mu(A) & \leqslant \frac{1}{\lambda}\Big( \int_{A} |T^{(1)}\mathbbm{1}_{E}|\mathrm{d}\mu + \int_{A} |T^{(2)}\mathbbm{1}_{E}|\mathrm{d}\mu\Big) \\
& \leqslant \frac{1}{\lambda}\Big( \|T^{(1)}\mathbbm{1}_{E}\|_{L^{q_{1}}}\mu(A)^{\frac{1}{q_{1}'}} + \|T^{(2)}\mathbbm{1}_{E}\|_{L^{q_{2}}}\mu(A)^{\frac{1}{q_{2}'}}\Big), \\
\end{aligned}
\end{equation} where, in keeping with the convention, $q_{1}'$ and $q_{2}'$ denote conjugates of $q_{1}$ and $q_{2}$ respectively.

For the finite sum $T^{(1)}$, we apply the first estimate \eqref{1}, the positive exponential bound, in the lemma, and for the remaining infinite $T^{(2)}$, we apply \eqref{2}, the negative exponential bound. These give \begin{equation} \begin{aligned}
\|T^{(1)}\mathbbm{1}_{E}\|_{L^{q_{1}}} 
& \leqslant \displaystyle\sum_{j=1}^{\rho} M_{1}2^{\beta_{1}j}\|\mathbbm{1}_{E}\|_{L^{p_{1}}}=\displaystyle\sum_{j=1}^{\rho} M_{1}2^{\beta_{1}j}\mu(E)^{\frac{1}{p_{1}}} \\
& \leqslant CM_{1}2^{\beta_{1}\rho}\mu(E)^{\frac{1}{p_{1}}}, \\
\end{aligned}
\end{equation}
and \begin{equation} \begin{aligned}
\|T^{(2)}\mathbbm{1}_{E}\|_{L^{q_{2}}} 
& \leqslant \displaystyle\sum_{j=\rho+1}^{\infty} M_{2}2^{-\beta_{2}j}\|\mathbbm{1}_{E}\|_{L^{p_{2}}}=\displaystyle\sum_{j=\rho+1}^{\infty} M_{2}2^{-\beta_{2}j}\mu(E)^{\frac{1}{p_{2}}} \\
& \leqslant CM_{2}2^{-\beta_{2}\rho}\mu(E)^{\frac{1}{p_{2}}}. \\
\end{aligned}
\end{equation}
Substituting the above two inequalities in \eqref{4}, we have \begin{equation}
\mu(A) \leqslant \frac{C}{\lambda}\Big(M_{1}2^{\beta_{1}\rho}\mu(E)^{\frac{1}{p_{1}}}\mu(A)^{\frac{1}{q_{1}'}}+M_{2}2^{-\beta_{2}\rho}\mu(E)^{\frac{1}{p_{2}}}\mu(A)^{\frac{1}{q_{2}'}}\Big).
\end{equation}
If we minimized the right-side of the last inequality by choosing appropriate $\rho$, then we would obtain \begin{equation} \label{5}
\mu(A) \leqslant \frac{C}{\lambda}M_{1}^{\theta}M_{2}^{1-\theta}\mu(E)^{\frac{1}{p}}\mu(A)^{1-\frac{1}{q}},
\end{equation} which would imply the estimate \eqref{6} that we set out to prove. The $\rho$ achieving the favorable result is such that \begin{equation} \label{7}
2^{\rho}=\Big( \frac{M_{1}}{M_{2}}\mu(E)^{\frac{1}{p_{2}}-\frac{1}{p_{1}}}\mu(A)^{\frac{1}{q_{2}'}-\frac{1}{q_{1}'}}\Big)^{\frac{1}{\beta_{1}+\beta_{2}}}.
\end{equation}

However, the $\rho$ here must be a nonnegative integer, so we need a few more lines to guarantee that the above reasoning goes through well. If \[
\frac{M_{1}}{M_{2}}\mu(E)^{\frac{1}{p_{2}}-\frac{1}{p_{1}}}\mu(A)^{\frac{1}{q_{2}'}-\frac{1}{q_{1}'}} > 1, \]
then we may let $\rho$ be the nonnegative integer such that \[
2^{\rho} < \Big( \frac{M_{1}}{M_{2}}\mu(E)^{\frac{1}{p_{2}}-\frac{1}{p_{1}}}\mu(A)^{\frac{1}{q_{2}'}-\frac{1}{q_{1}'}}\Big)^{\frac{1}{\beta_{1}+\beta_{2}}} \leqslant 2^{\rho+1}. \] This choice of $\rho$ will lead us to the same estimate as \eqref{5}, with a worse $C$. If however \[
\frac{M_{1}}{M_{2}}\mu(E)^{\frac{1}{p_{2}}-\frac{1}{p_{1}}}\mu(A)^{\frac{1}{q_{2}'}-\frac{1}{q_{1}'}} \leqslant 1, \]
then we may simply let $\rho$ be $0$ and sum the negative exponential bounds. The result is \begin{equation}
\mu(A) \leqslant \frac{1}{\lambda}M_{2}\mu(E)^{\frac{1}{p_{2}}}\mu(A)^{\frac{1}{q_{2}'}}.
\end{equation}
Rearranging terms, this is equivalent to \[
\lambda\mu(A)^{\frac{1}{q}} \leqslant M_{1}^{\theta}M_{2}^{1-\theta}\mu(E)^{\frac{1}{p}} \cdot \Big( \frac{M_{1}}{M_{2}}\mu(E)^{\frac{1}{p_{2}}-\frac{1}{p_{1}}}\mu(A)^{\frac{1}{q_{2}'}-\frac{1}{q_{1}'}}\Big)^{\frac{\beta_{2}}{\beta_{1}+\beta_{2}}}. \]
Considering our assumption, the above is in fact a stronger result than we need. We thence concludes the proof of Bourgain's interpolation technique.

\vspace{0.5cm}

Now we embark on the task of proving Proposition 1. The essence in the proof is the following lemma on an oscillatory integral, which is an extension of Proposition 2.2 in Huang-Sogge \cite{HS}.

\begin{lemma}
Suppose $g$ defines a smooth Riemannian metric on $\mathbb{R}^{n}$ that is close to the Euclidean one and has injectivity radius larger than 10. If the function \[a(x, y)\in C^{\infty} \big( B_{2}(O) \times B_{2}(O) \setminus \{(x, x): x \in B_{2}(O) \} \big)\]
satisfies the estimates \begin{equation} \label{8}
a(x, y) \leqslant Cd_{g}(x, y)^{2-n}, \quad \text{if} \quad d_{g}(x, y) \leqslant \frac{1}{\lambda},
\end{equation}
\begin{equation} \label{9}
\big|\partial^{\alpha}_{x, y}a(x, y)\big| \leqslant C_{\alpha}\lambda^{\frac{n-3}{2}}d_{g}(x, y)^{-\frac{n-1}{2}-|\alpha|}, \quad \text{if} \quad d_{g}(x, y) \geqslant \frac{1}{\lambda},
\end{equation}
then, for exponents as in the Proposition 1, i.e. \[
\frac{1}{r}-\frac{1}{s}=\sigma, \quad \frac{2}{n+1} \leqslant \sigma \leqslant \frac{2}{n}, \]
\[\frac{2n}{n-1+2n\sigma} < r < \frac{2n}{n+1}, \]
we have \begin{equation} \label{10}
\Big\| \int e^{i \lambda d_{g}(x, y)}a(x, y)f(y)\mathrm{d}y \Big\|_{L^{s}(B_{1}(O))} \leqslant C\lambda^{n\sigma-2}\|f\|_{L^{r}(B_{1}(O))},
\end{equation}
where $f \in C^{\infty}_{0}(\mathbb{R}^{n})$ is supported in $B_{1}(O)$.
\end{lemma}

In the lemma, $B_{2}(O)$ and $B_{1}(O)$ are balls centered at the origin with respect to the metric $g$, and $d_{g}(x, y)$ is the Riemannian distance function induced by $g$.

\paragraph{Proof of Lemma 2}

Recall that for a fixed line segment \[
\frac{1}{r}-\frac{1}{s}=\sigma, \quad \frac{2}{n+1} \leqslant \sigma \leqslant \frac{2}{n}, \]
\[\frac{2n}{n-1+2n\sigma} < r < \frac{2n}{n+1}, \]
we aim to prove restricted weak type inequality at the endpoints using Bourgain's interpolation technique, and then apply real interpolation to obtain strong estimate for exponents in between. By duality, it suffices to deal with one endpoint, say \[(\frac{1}{p}, \frac{1}{q})=(\frac{n+1}{2n}, \frac{n+1-2n\sigma}{2n}). \] In what follows, $(p, q)$ will specifically refer to the above endpoint. 

We dyadically decompose the kernel $a(x, y)$ of the oscillatory integral in the lemma. So choose a Littlewood-Paley bump function $\beta(t) \in C_{0}^{\infty}\big( (\frac{1}{2}, 1) \big)$ such that \[ \displaystyle\sum_{j=-\infty}^{\infty} \beta(2^{-j}t)=1 \] 
whenever $t > 0$. Denote the operator given by the oscillatory integral in the lemma as $T$, and define operators $T_{j}$, $j=0, 1, 2, \cdots$, by
\[ T_{j}f(x)=\int e^{i \lambda d_{g}(x, y)}\beta(\lambda 2^{-j}d_{g}(x, y))a(x, y)f(y) \mathrm{d}y, \quad k=1, 2, \cdots, \]
\[ T_{0}=T-\displaystyle\sum_{j=1}^{\infty} T_{j}. \]
$T_{0}$ can be handled by Young's inequality, if we notice estimate \eqref{8} on $a(x, y)$ and the support property of $\beta(t)$: \begin{equation} \label{11} \begin{aligned}
\|T_{0}f(x)\|_{L^{q}(B_{1}(O))}
& \leqslant C\|f\|_{L^{p}(B_{1}(O))} \Big( \int_{d_{g}(x, y) \leqslant \frac{1}{\lambda}} \big( d_{g}(x, y)^{2-n} \big) ^{\frac{1}{1-\sigma}} \mathrm{d}y \Big) ^{1-\sigma} \\
& = C\lambda^{n\sigma-2}\|f\|_{L^{p}(B_{1}(O))}. \\
\end{aligned}
\end{equation}
Note that the condition $\sigma \leqslant \frac{2}{n}$ is necessary to ensure a finite bound in this estimate. 

For $T_{k}$, we apply Bourgain's result Lemma 1. Recall the line segment $s=\frac{n+1}{n-1}r'$, $1 \leqslant r \leqslant 2$ on which Stein's oscillatory integral theorem (see \cite{Stein}) holds. This line segment intersects with the one we are interested in, \[
\frac{1}{r}-\frac{1}{s}=\sigma, \quad \frac{2n}{n-1+2n\sigma} < r < \frac{2n}{n+1}, \]
at \[
P: \Big( \frac{n+1}{2n}\sigma+\frac{n-1}{2n}, -\frac{n-1}{2n}\sigma+\frac{n-1}{2n} \Big).\]
Since the oscillatory integral in Lemma 2 is closely related to the oscillatory integral in Stein's theorem, we are naturally led to considering interpolating between the above point $P$ of intersection with the point $Q(\sigma, 0)$, also on the line $\frac{1}{r}-\frac{1}{s}=\sigma$ but at the same time lying on the horizontal axis. These two points are shown in Figure \ref{fig2} below. For simplicity of notation, we denote the exponents corresponding to $P$ by $(p_{1}, q_{1})$ and those corresponding to $Q$ by $(p_{2}, q_{2})$, so that $P=(\frac{1}{p_{1}}, \frac{1}{q_{1}})$ and $Q=(\frac{1}{p_{2}}, \frac{1}{q_{2}})$.  Note that at $Q$, we can always obtain a trivial bound from Holder's inequality.

There is however, a small obstacle that need be overcome. When $\sigma=\frac{2}{n+1}$, the point $P$ coincides with the endpoint at which we want to prove restricted weak type estimate, hence in this case, interpolating in the above way does not help. This is fortunately, not too much of a trouble, and we are able to remedy it by interpolating instead between the two endpoints $A$ and $B$ of the line segment in Stein's oscillatory integral theorem. See Figure \ref{fig2} below. In the rest of the proof therefore, we first treat the general case, and then turn to the exceptional one.

\begin{figure}
  \centering
    \includegraphics[height=6.5cm]{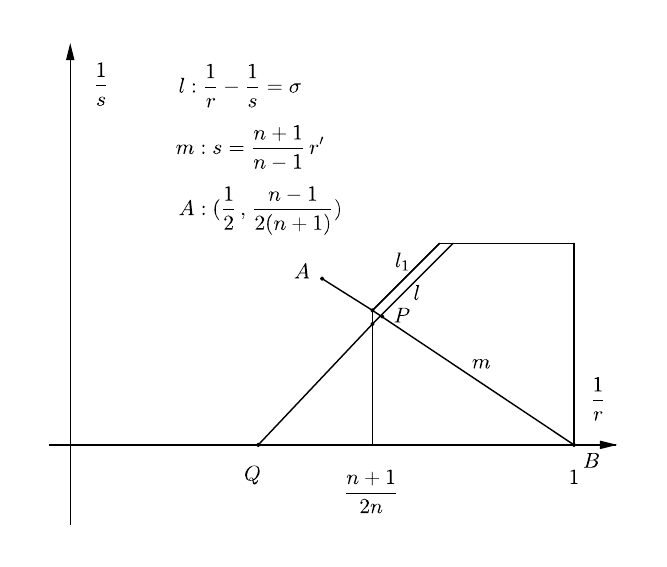}
\caption{Interpolation argument 1}
  \label{fig2}
\end{figure}  

Now we begin the interpolation process. At $P$, our goal is to obtain the following bound for every $j$: \begin{equation} \label{12}
\|T_{j}f\|_{L^{q_{1}}(B_{1}(O))} \leqslant C2^{\big( -\frac{(n+1)\sigma}{2}+1\big)j}\lambda^{n\sigma-2}\|f\|_{L^{p_{1}}(B_{1}(O))},
\end{equation}
where the bound $C$ is independent of $j$. By a dilation argument, denoting $\epsilon=\frac{2^{j}}{\lambda}$, this would be a consequence of the estimate \begin{equation} \label{13}
\|S_{j}f\|_{L^{q_{1}}(B_{\epsilon}(O))} \leqslant C2^{\big( \frac{n-1}{2}\sigma-\frac{n-1}{2}\big) j}\|f\|_{L^{p_{1}}(B_{\epsilon}(O))}
\end{equation} 
on the operator $S_{j}$ defined by \begin{equation} \label{24}
S_{j}f(x)=\int e^{i 2^{j}\phi_{j}(x, y)}b_{j}(x, y)f(y)\mathrm{d}y,
\end{equation}
\[\phi_{j}(x, y)=\frac{1}{\epsilon}d_{g}(\epsilon x, \epsilon y),\]
\[b_{j}(x, y)=2^{\frac{n-1}{2}j}\lambda^{2-n}\beta(\phi_{j}(x, y))a(\epsilon x, \epsilon y). \]
Notice that if we introduce the new metric $g(\epsilon x)$ by dilation, then $\phi_{j}(x, y)$ is precisely the  Riemannian distance function this new metric induces. Also, the ball of radius $\epsilon$ have radius $1$ under $g(\epsilon x)$. Furthermore, the amplitude $b_{j}(x, y)$ vanishes if $\phi_{j}(x, y) \notin [\frac{1}{2}, 1]$, and satisfies the estimate \[
\big| \partial^{\alpha}_{x, y}b_{j}(x, y)\big| \leqslant C_{\alpha}, \]
where the bound $C_{\alpha}$ is independent of $j$. These observations tell that the oscillatory integral defined by the operator $S_{j}$ satisfies the hypotheses of Stein's oscillatory integral theorem. Then \eqref{13}, the equation we wish to prove, follows directly from Stein. We point out here that the condition $\sigma < \frac{2}{n+1}$ is necessary to ensure that the power $-\frac{(n+1)\sigma}{2}+1$ on $2$ is negative. The proof of the estimate \eqref{12} at point $P$ has been accomplished.

As mentioned above, at $Q$, there is the trivial bound resulting from Holder's inequality. Noting the condition \eqref{9} on $a(x, y)$ when $d_{g}(x, y)>\frac{1}{\lambda}$, we readily get \begin{equation} \label{14}
\|T_{j}f\|_{L^{q_{2}}(B_{1}(O))} \leqslant C2^{\big(\frac{n+1}{2}-n\sigma\big)j}\lambda^{n\sigma-2}\|f\|_{L^{p_{2}}(B_{1}(O))}.
\end{equation}
When $\frac{2}{n+1} < \sigma \leqslant \frac{2}{n}$, the power $\frac{n+1}{2}-n\sigma>0$.

It remains only some elementary computation.
\[ \theta=\frac{\frac{n+1}{2}-n\sigma}{\frac{n+1}{2}-n\sigma+\frac{(n+1)\sigma}{2}-1}, \]
\[\frac{\theta}{p_{1}}+\frac{1-\theta}{p_{2}}=\frac{n+1}{2n}=\frac{1}{p}, \] 
\[\frac{\theta}{q_{1}}+\frac{1-\theta}{q_{2}}=\frac{n+1}{2n}-\sigma=\frac{1}{q}. \]
The computation shows that we have arrived exactly at $(\frac{1}{p}, \frac{1}{q})$, the endpoint that is our goal. Thus, Bourgain's interpolation technique gives us the desired restricted weak type inequality for the sum of $T_{j}$. Together with the strong estimate \eqref{11} for $T_{0}$, this yields \[
\|Tf\|_{L^{q, \infty}(B_{1}(O))} \leqslant C\lambda^{n\sigma-2}\|f\|_{L^{p, 1}(B_{1}(O))}. \]
Duality gives the same inequality for the other endpoint. Finally, real interpolation produces the conclusion in the lemma.

There is still the exceptional case $\frac{1}{r}-\frac{1}{s}=\frac{2}{n+1}$ that need be tackled. As mentioned before, we interpolate instead between the two endpoints $A(\frac{1}{2}, \frac{n-1}{2(n+1)})$ and $B(1, 0)$ in Stein's oscillatory integral theorem. The procedure is pretty much the same as with the general case, so we provide an outline. For convenience, we still denote the exponents corresponding to $A$ and $B$ as $(p_{1}, q_{1})$ and $(p_{2}, q_{2})$, respectively. At $A$, we wish to show for each $j$, \begin{equation} \label{25}
\|T_{j}f\|_{L^{q_{1}}(B_{1}(O))} \leqslant C2^{\frac{1}{2}j}\lambda^{-\frac{n+2}{n+1}}\|f\|_{L^{p_{1}}(B_{1}(O))},
\end{equation} where the $C$ is independent of $j$. By a dilation argument, this follows from \begin{equation} \label{27}
\|S_{j}f\|_{L^{q_{1}}(B_{1}(O))} \leqslant C2^{-\frac{n^{2}-n}{n+1}}\|f\|_{L^{p_{1}}(B_{1}(O))},
\end{equation}
where $S_{j}$ is as in \eqref{24}, and the unit ball $B_{1}(O)$ is with respect to the dilated metric $g(\epsilon x)$. But \eqref{27} is a direct consequence of Stein's oscillatory integral theorem, so our work at point $A$ is done. At $B$, we have the trivial estimate \begin{equation} \label{28}
\|T_{j}f\|_{L^{q_{2}}(B_{1}(O))} \leqslant C2^{-\frac{n-1}{2}j}\lambda^{n-2}\|f\|_{L^{p_{2}}(B_{1}(O))}.
\end{equation}
Note again that the power $-\frac{n-1}{2}$ of $2$ in \eqref{28} is negative. Then we verify \[
\theta=\frac{\frac{n-1}{2}}{\frac{n-1}{2}+\frac{1}{2}},\]
\[\theta (\frac{1}{p_{1}}, \frac{1}{q_{1}})+(1-\theta)(\frac{1}{p_{2}}, \frac{1}{q_{2}})=(\frac{n+1}{2n}, \frac{(n-1)^{2}}{2n(n+1)}). \]
Once more, this latter pair of exponents correspond exactly to one endpoint of the extraordinary line segment. Therefore, restricted weak type inequality holds at this endpoint \begin{equation} \label{26}
\|Tf\|_{L^{\frac{2n(n+1)}{(n-1)^{2}}, \infty}(B_{1}(O))} \leqslant C\lambda^{-\frac{2}{n+1}}\|f\|_{L^{\frac{2n}{n+1}, 1}(B_{1}(O))}.
\end{equation}
The rest of the work, i.e. an application of duality and real interpolation, is the same as in the general case. This finishes the entire proof of Lemma 2.

\paragraph{Proof of Proposition 1}

With Lemma 2, we can proceed to show Proposition 1. The procedure parallels that in Huang-Sogge \cite{HS}. By the compactness of the sphere, we may assume all functions $f$ appearing in this proof to be supported in a ball of radius $1$. It is well-known that $H_{k}(x, y) \leqslant Ck^{n-1}$ (\cite{Sogge}, $\S$ 3.4). However, to prove the proposition, we need a more precise bound on $H_{k}(x, y)$, at least for large $k$. The result we cited is Proposition 2.1 in Huang-Sogge \cite{HS}.

\begin{proposition}
When $k$ is large, \begin{multline} \label{15}
H_{k}(x, y)=\lambda_{k}^{\frac{n-1}{2}}\displaystyle\sum_{\pm} a_{\pm}(k; x, y)e^{\pm i\lambda_{k}d_{S^{n}}(x, y)}, \\
\text{for} \quad \frac{1}{\lambda_{k}} \leqslant d_{S^{n}}(x, y) \leqslant \frac{3\pi}{4},
\end{multline}
where $a_{\pm}(k; x, y)$ are smooth functions satisfying for each $j=0, 1, 2, \cdots$, \begin{equation}
\big|\partial^{j}_{x, y}a_{\pm}(k; x, y)\big| \leqslant Cd_{S^{n}}(x, y)^{-j}, \quad \text{if} \quad \frac{1}{\lambda_{k}} \leqslant d_{S^{n}}(x, y) \leqslant \frac{3\pi}{4}.
\end{equation}
We also have the expression \begin{multline} \label{16}
H_{k}(x, y)=(-1)^{k}\lambda_{k}^{\frac{n-1}{2}}\displaystyle\sum_{\pm} a_{\pm}(k; x, y^{\ast})e^{\pm i\lambda_{k}d_{S^{n}}(x, y^{\ast})},\\
\text{for} \quad \frac{\pi}{4} \leqslant d_{S^{n}}(x, y) \leqslant \pi-\frac{1}{\lambda_{k}}.
\end{multline}
\end{proposition}

Choose $\alpha(t) \in C^{\infty}(\mathbb{R}_{+})$ such that \[
\alpha(t)=1 \quad \text{if} \quad t \leqslant \delta; \quad \alpha(t)=0 \quad \text{if} \quad t \geqslant 2\delta, \]
where $\delta$ is to be specified later. Define \[
\tilde{H}_{k}(x, y)=\alpha(d_{S^{n}}(x, y))H_{k}(x, y).\]
If $\delta$ is sufficiently small, then by Proposition 2 and the estimate $H_{k}(x, y) \leqslant Ck^{n-1}$, it is easy to see that $\frac{\tilde{H}_{k}(x, y)}{\lambda_{k}}$ satisfies that hypotheses set for the amplitude $a(x, y)$ in Lemma 2. Therefore, applying that lemma yields \begin{equation} \label{17}
\|\tilde{H}_{k}f\|_{L^{s}(S^{n})} \leqslant k^{n\sigma-1}\|f\|_{L^{r}(S^{n})},
\end{equation}
where the exponents $r$ and $s$ are on our fixed line segment \eqref{30}. For the same reason, if we define \[
\tilde{H}_{k}^{\ast}(x, y)=\alpha(d_{S^{n}}(x, y^{\ast}))H_{k}(x, y), \]
we have \begin{equation} \label{18}
\|\tilde{H}_{k}^{\ast}f\|_{L^{s}(S^{n})} \leqslant k^{n\sigma-1}\|f\|_{L^{r}(S^{n})}.
\end{equation}

What remains is $U_{k}=H_{k}-\tilde{H}_{k}-\tilde{H}_{k}^{\ast}$, an operator given by the oscillatory integral \[
U_{k}f(x)=\lambda_{k}^{\frac{n-1}{2}}\displaystyle\sum_{\pm} \int_{S^{n}} \alpha_{\pm}(k; x, y)e^{\pm i\lambda_{k}d_{S^{n}}(x, y)}f(y) \mathrm{d}y, \]
where the functions $\alpha_{\pm}(k; x, y)$ vanishes when $d_{S^{n}}(x, y) \notin [\delta, \pi-\delta]$, and satisfies \[
|\partial^{\gamma}_{x, y}\alpha_{\pm}(k; x, y)| \leqslant C_{\gamma}, \]
with the bound $C_{\gamma}$ independent of $k$. This reminds us to utilize again Stein's oscillatory integral theorem. We apply it to the pair of exponents $\big( \frac{n+1}{2n}, \frac{(n-1)^{2}}{2n(n+1)}\big)$,
which is point $C$ in Figure \ref{fig3} below. The result is \begin{equation} \label{19}
\|U_{k}f\|_{L^{\frac{2n(n+1)}{(n-1)^{2}}}(S^{n})} \leqslant C\lambda_{k}^{\frac{n-1}{n+1}}\|f\|_{L^{\frac{2n}{n+1}}(S^{n})}.
\end{equation}
Furthermore, at $D(\frac{n+1}{2n}, 0)$, the point on the horizontal axis also labeled in Figure \ref{fig3}, we have an estimate from Holder's inequality \begin{equation}
\|U_{k}f\|_{L^{\infty}(S^{n})} \leqslant C\lambda_{k}^{\frac{n-1}{2}}\|f\|_{L^{\frac{2n}{n-1}}(S^{n})}.
\end{equation}
Interpolating between $C$ and $D$ then gives the desired strong estimate at one endpoint of the line segment \eqref{30}: \begin{equation} \label{20}
\|U_{k}f\|_{L^{\frac{2n}{n+1-2n\sigma}}(S^{n})} \leqslant C\lambda^{n\sigma-1}\|f\|_{L^{\frac{2n}{n+1}}(S^{n})}.
\end{equation} As before, duality produces the inequality at the other endpoint, and another interpolation results in the estimate we are seeking for exponents between the endpoints. Combining \eqref{17}, \eqref{18} and \eqref{20}, we finish our proof of Proposition 1.

\begin{figure}
  \centering
    \includegraphics[height=6.5cm]{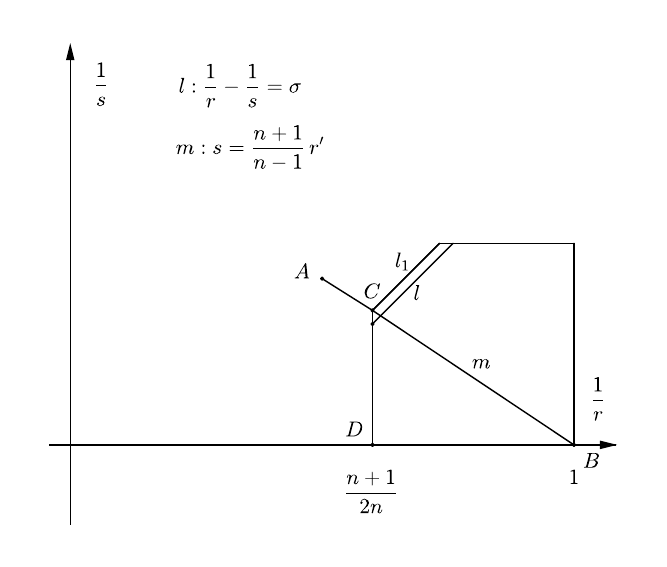}
\caption{Interpolation argument 2}
  \label{fig3}
\end{figure}

\maketitle
\section{Proof of Theorem 1}
We finally come to the proof of the main theorem of this article. It follows the lines in Huang-Sogge \cite{HS}. Let $P$ represent the operator $\sqrt{-\Delta_{S^{n}}+(\frac{n-1}{2})^{2}}$. We apply the formula \begin{equation} \label{21}
(-P^{2}+(\lambda+i\mu)^{2})^{-1}=\frac{\mathrm{sgn}\mu}{i(\lambda+i\mu)}\int _{0}^{\infty} e^{i(\mathrm{sgn}\mu)\lambda t}e^{-|\mu|t}\mathrm{cos}(tP) \mathrm{d}t,
\end{equation} 
which can be found in \cite{BSSY}. 

Choose a nonnegative even function $\rho(t) \in C_{0}^{\infty}(\mathbb{R})$ that satisfies \[
\rho(t)=1 \quad \text{if} \quad |t|<\frac{1}{2}; \quad \rho(t)=0 \quad \text{if} \quad |t|\geqslant 1. \] By Proposition 2.4 in Huang-Sogge \cite{HS}, the kernel of the operator \[
R_{0}f=\frac{\mathrm{sgn}\mu}{i(\lambda+i\mu)}\int _{0}^{\infty} \rho(t)e^{i(\mathrm{sgn}\mu)\lambda t}e^{-|\mu|t}\mathrm{cos}(tP) \mathrm{d}t \]
has the expression \[
\displaystyle\sum_{\pm} b_{\pm}(\lambda; x, y)e^{\pm i\lambda d_{S^{n}}(x, y)}+O\big((d_{S^{n}}(x, y))^{2-n}\big), \]
where the functions $b_{\pm}(\lambda; x, y)$ vanish when $d_{S^{n}}(x, y)$ is near to $\pi$, and satisfy the estimates \eqref{8} and \eqref{9} in Lemma 2, with all constants independent of $\lambda \geqslant 1$. Therefore, Lemma 2 and Young's inequality give \begin{equation} \label{22}
\|R_{0}f\|_{L^{s}(S^{n})} \leqslant C_{r, s}\lambda^{n\sigma-2}\|f\|_{L^{r}(S^{n})}.
\end{equation}

For the remaining ``$1-\rho(t)$'' part, denote \[
m_{\lambda, \mu}(\tau)=\frac{\mathrm{sgn}\mu}{i(\lambda+i\mu)}\int _{0}^{\infty} (1-\rho(t))e^{i(\mathrm{sgn}\mu)\lambda t}e^{-|\mu|t}\mathrm{cos}(t\tau) \mathrm{d}t. \]
Then, an easy integration by parts argument shows \[
|m_{\lambda, \mu}(\tau)| \leqslant C_{N}\lambda^{-1}(1+|\lambda-\tau|)^{-N} \]
for every natural $N$, whenever $\tau \geqslant 0$, and $\lambda, |\mu| \geqslant 1$. Therefore, by the spectral theorem, \begin{equation} \label{23} \begin{aligned}
\Big\| \big( (-P^{2}+(\lambda+i\mu)^{2})^{-1}-R_{0}\big) f & \Big\|_{L^{s}(S^{n})} =\Big\| \displaystyle\sum _{k=0}^{\infty} m_{\lambda, \mu}(\lambda_{k})H_{k}f\Big\|_{L^{s}(S^{n})} \\
& \leqslant \displaystyle\sum _{k=0}^{\infty} \|m_{\lambda, \mu}\|_{L^{\infty}}\|H_{k}f\|_{L^{s}(S^{n})} \\
& \leqslant C \displaystyle\sum_{k=0}^{\infty} \lambda^{-1}(1+|\lambda-k|)^{-3}k^{n\sigma-1}\|f\|_{L^{r}(S^{n})} \\
& \leqslant C\lambda^{n\sigma-2}\|f\|_{L^{r}(S^{n})}.\\
\end{aligned}
\end{equation}
This concludes the proof of Theorem 1.

\nocite{Gutierrez}

\vspace{0.5cm}

\textbf{Acknowledgement}. The author would like to express his gratitude to his advisor, Professor Christopher D. Sogge, for the invaluable guidance he provided and the illuminating discussions with him.

\bibliography{Manifold}

\bibliographystyle{plain}

\vspace{0.5cm}

\setlength{\parindent}{0cm}
\textit{Email address}: tyren@math.jhu.edu

\setlength{\parindent}{0cm}
DEPARTMENT OF MATHEMATICS, JOHNS HOPKINS UNIVERSITY, BALTIMORE, MD 21218, USA

\end{document}